\documentstyle[epsfig]{article}

\newtheorem{theorem}{Theorem}
\newtheorem{lemma}[theorem]{Lemma}
\newtheorem{prop}[theorem]{Proposition}
\newtheorem{cor}[theorem]{Corollary}

\title{Line transversals to disjoint balls}
\author{Ciprian Borcea\thanks{Rider University, Lawrenceville, NJ 08648, \texttt{borcea@rider.edu}.}
\and Xavier Goaoc\thanks{LORIA - INRIA Lorraine, Nancy, France, \texttt{goaoc@loria.fr}.} \and Sylvain
Petitjean\thanks{LORIA - CNRS, Nancy, France, \texttt{petitjea@loria.fr}.}}
\date{}

\begin{document}
\maketitle

\begin{abstract}
We prove that the set of directions of lines intersecting three disjoint balls in $R^3$ in a given order is a
strictly convex subset of $S^2$. We then generalize this result to $n$ disjoint balls in $R^d$. As a
consequence, we can improve upon several old and new results on line transversals to disjoint balls in arbitrary
dimension, such as bounds on the number of connected components and Helly-type theorems.
\end{abstract}

\section{Introduction}

\medskip \noindent
Helly's theorem \cite{H} of 1923 opened a large field of inquiry now designated as {\em geometric transversal
theory}. A typical concern is the study of all $k$-planes (also called $k$-flats) which intersect all sets of a
given family of subsets (or objects) in $R^d$. These are the $k$-transversals of the given family and they
define a certain subspace of the corresponding Grassmannian. True to its origin, transversal theory usually
implicates {\em convexity} in some form: either in its assumptions, or in its proofs, or, most likely, in both.

\medskip \noindent
In what follows, $k=1$ and the objects will be disjoint closed balls with arbitrary radii in $R^d$. Our main
result is the following convexity theorem:

\begin{theorem}\label{thm:convexity}
The directions of all oriented lines intersecting a given finite family of disjoint balls in $R^d$ in a specific
order form a strictly convex subset of the sphere $S^{d-1}$.
\end{theorem}

\noindent
with the immediate consequence that the connected components in the space of line transversals correspond with
all possible {\em geometric permutations} of the given family, where a geometric permutation is understood as a
pair of orderings defined by a single line transversal with its two orientations.

\newpage
\noindent Before discussing other implications, we want to emphasize that the {\em key} to our theorem
resides in the case of {\em three disjoint balls in} $R^3$, and the approach we use to settle this case is
geometrically quite revealing, in that it shows the nuanced dependency of the convexity property on the {\em
curve of common tangents} to the three bounding spheres.

\subsection{Relations with previous work}

\medskip \noindent
Helly's theorem \cite{H} states that a finite family ${\cal S}$ of convex sets in $R^d$ has non-empty
intersection if and only if any subfamily of size at most $d+1$ has non-empty intersection. Passing from $k=0$
to $k=1$, one of the early results is due to Danzer \cite{D} who proved that $n$ disjoint {\bf unit} disks in
the plane have a line transversal if and only if every five of them have a line transversal. Hadwiger's theorem
\cite{Had1}, which allows arbitrary disjoint convex sets in the plane as objects, showed the importance of the
{\em order} in which oriented line transversals meet the objects: when every three objects have an oriented line
transversal respecting some fixed order of the whole family, there must be a line transversal for the family.

\medskip \noindent
This stimulated the interest of comparing, in arbitrary dimension, two equivalence relations for line
transversals: the coarser one, {\em geometric permutation}, determined by the order in which the given disjoint
objects are met (up to reversal of orientation), and the finer one, {\em isotopy}, determined by the connected
components of the space of transversals.

\medskip \noindent
In general, for $d\geq 3$, the ``gap'' between the two notions may be wide \cite{GPW}, and families for which
the two notions coincide are thereby ``remarkable''. The first examples of such families are ``thinly
distributed'' balls\footnote{A family of balls is {\em thinly distributed} if the distance between the centers
of any two balls is at least twice the sum of their radii.} in arbitrary dimension, as observed by Hadwiger
\cite{Had2}. Then, the work of Holmsen et al. \cite{HKL} showed that disjoint {\bf unit} balls in $R^3$ provide
``remarkable'' cases as well. They verified the convexity property in the case of equal radii, and their method
can be extended to the larger class of ``pairwise inflatable" balls\footnote{A family of balls is {\em pairwise
inflatable} if the squared distance between the centers of any two balls is at least twice the sum of their
squared radii.} in arbitrary dimension \cite{CGHP}, inviting the obvious question regarding disjoint balls of
arbitrary radii. The significance of this problem is also discussed in the recent notes \cite[pg~191--195]{PS}
where one can find ampler references to related literature.

\medskip \noindent
Our solution for the case of arbitrary radii is based on a new approach, suggested by the detailed study of the
curve of common tangents to three spheres in $R^3$~\cite{Bor2}. The main ideas are outlined in
Section~\ref{sec:outline} as a preamble to the detailed proof in Section~\ref{sec:computation}.

\medskip \noindent
In dimension three, particularly, there are  connections with other problems in visibility and geometric
computing. Changes of visibility (or ``visual events'') in a scene made of smooth obstacles typically occur for
multiple tangencies between a line and some of the obstacles \cite{Pl}. Tritangent and quadritangent lines play
a prominent role in this picture, as they determine the 1- and 0-dimensional faces of visibility structures. An
attractive case is that of four balls in $R^3$, which allow, generically, up to twelve common real tangents
\cite{MPT}. Degenerate configurations are identified in \cite{BGLP}.  Variations on such problems, where
reliance on algebraic geometry comes to the forefront, are surveyed in \cite{ST}. See also a brief account in
\cite{Bor1}.

\subsection{Further implications}

\medskip \noindent
Danzer's theorem \cite{D} motivated several other attempts to generalize Helly's result for $k=1$, i.e. for line
transversals. Whereas Helly's theorem only requires convexity, the case $k=1$ appears to be more sensitive to
the geometry of the objects. In particular, Holmsen and Matou{\v s}ek \cite{HM} showed that no such theorem
holds in general for families of disjoint translates of a convex set -- not even with restriction on the
ordering \emph{\`a la Hadwiger}. Our Theorem~\ref{thm:convexity} has consequences in this direction, presented
below in Section~\ref{sec:implications}.

\medskip \noindent
Hadwiger's proof of his Transversal Theorem \cite{Had1} relies on the observation that any \emph{minimal pinning
configuration}, i.e. family of objects with an isolated line transversal that would become non-isolated should
any of the objects be removed, has size $3$ if the objects are disjoint convex sets in the plane.
Theorem~\ref{thm:convexity} implies that any minimal pinning configuration of disjoint balls in $R^d$ has size
at most $2d-1$ (Corollary~\ref{cor:pinning}). A generalization of Hadwiger's theorem for families of disjoint
balls then follows (Corollary~\ref{cor:Hadwiger}).

\section{Preliminaries}\label{sec:preliminaries}

\paragraph{Notations and prerequisites.}
For any two vectors $a$, $b$ of $R^3$, we denote by $<a,b>$ their dot product and by $a \times b$ their cross
product.

\medskip \noindent
The space of directions in $R^3$ is the real projective space $P_2=P_2(R)$ envisaged either as the space of
lines through the origin (and then the direction of a line is given by its parallel through the origin), or as
the ``hyperplane at infinity'' in the completion $P_3=R^3\sqcup P_2$ (and then the direction of a line is simply
its point of intersection with the hyperplane at infinity). Convexity in $P_2$ is relative to the metric induced
from the standard metric of the sphere through the identification $S^2/Z_2=P_2$. All considerations can be
pulled-back to $S^2$ by orienting the lines.

\medskip \noindent
In following our convexity arguments related to three disjoint balls in $R^3$, it may be helpful to bear in mind
that the regions of $P_2$ determined by directions of line transversals are always contained in the
simply-connected side of some smooth conic (which is homeomorphic with a disc, while the other side is
homeomorphic with a M\"{o}bius band). When testing convexity, one may use affine charts $R^2$, and verify
locally, then globally, that the boundary curve ``stays on the same side of its tangent''. If this property were
to fail at some point, one must have an {\em inflection point} there, or, in one word, a {\bf flex}.

\medskip \noindent
We denote by $B_0, B_1, B_2$ three balls in $R^3$ with respective centers $c_0, c_1, c_2$ and squared radii
$s_0, s_1, s_2$, $s_k = r_k^2$. Since the degenerate case of collinear centers is easily obtained from the
generic case, we assume that we have a non-degenerate {\em triangle of centers}.

\paragraph{Direction-sextic.}
The directions of common tangent lines to $B_0,B_1, B_2$ make up an algebraic curve of degree six in $P_2$,
which we call the \emph{direction-sextic} and denote $\sigma$. To take advantage of symmetries in expressing
$\sigma$, we introduce the edge vectors $e_{ij}=c_j-c_i$, denote by $\delta_{ij} = <e_{ij},e_{ij}>$ their
squared norms and put

$$  q=q(u)=<u, u>,\\ \ \ \ \
  t_{ij}=t_{ji}=<e_{ij} \times u, e_{ij} \times u> = \delta_{ij}q-<e_{ij},u>^2 $$

\noindent Thus in $P_2(C)$, the equation $t_{ij}=0$ gives the two tangents from $e_{ij}$ to the imaginary conic
$q=0$.

\begin{prop}\label{sextic}
The direction-sextic for $B_0,B_1, B_2$ can be given by means of the Cayley determinant:

$$     \sigma=\sigma(u)= \det{\left( \begin{array}{ccccc}
           0 & 1 & 1 & 1 & 1 \\
           1 & 0    & qs_0   & qs_1   & qs_2 \\
           1 & qs_0 & 0      & t_{01} & t_{02} \\
           1 & qs_1 & t_{01} & 0      & t_{12} \\
           1 & qs_2 & t_{02} & t_{12} & 0
         \end{array} \right)} = 0 $$
\end{prop}

\medskip \noindent
{\em Proof:}\ One way to find the equation of the direction curve is to begin with a description of lines in
$R^3$ by parameters $(p,u)\in R^3\times P_2$, where $p$ is the orthogonal projection of the origin on the given
line, and $u$ is the direction of the line.  With $c_0=0$ and abbreviations:

$$ a_i=a_i(u)=<c_i\times u,c_i\times u>+(s_0-s_i)<u,u>=t_{0i}+(s_0-s_i)q, \ \
 i=1,2  $$

\noindent affine common tangents obey the system (see e.g. \cite{BGLP} or \cite{MPT}):

$$   <p,c_i>=\frac{a_i(u)}{2<u,u>}, \quad i=1,2,\ \
   <p,u>=0,\ \
   <p,p>=s_0 $$

\medskip \noindent
The direction-sextic is obtained by eliminating $p$ from this system. The fact that the resulting equation
allows the stated Cayley determinant expression is given a natural explanation in \cite{Bor2}, but can be
directly verified by computation. \ \ $\Box$

\medskip \noindent
The direction of an {\em oriented line} can be represented either by a point on the unit sphere, or by the whole
{\em ray} emanating from the origin and passing through that point. Our expression ``cone of directions'' stems
from the latter representation, which converts questions of convexity in $S^2$ into equivalent questions of
convexity in $R^3$. In the projective context, it will be understood that we mean the image via $S^2/Z_2=P_2$.

\begin{figure}[!t]
  \centerline{
    $\vcenter{\hbox{\epsfig{file=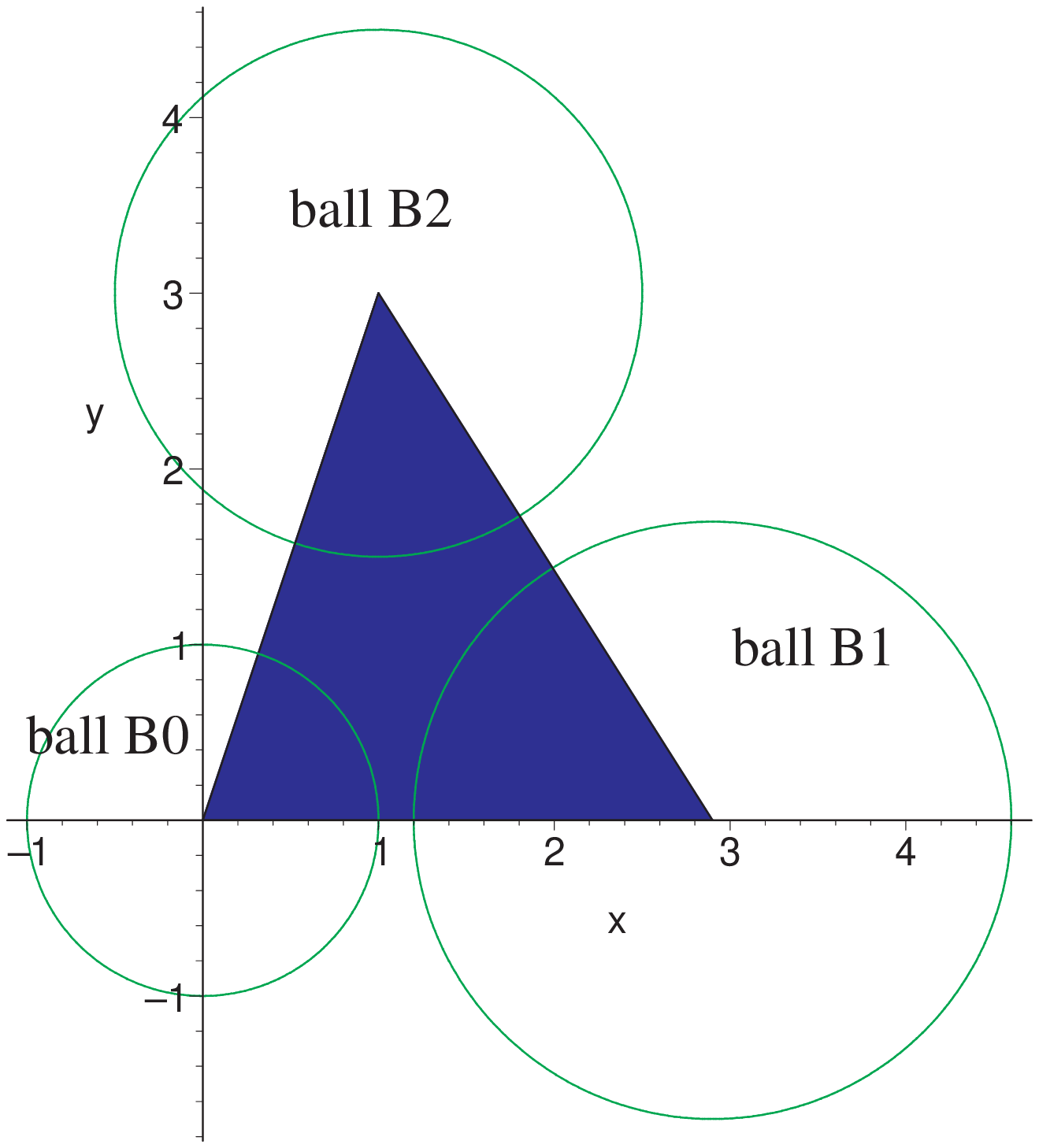,height=6cm,clip=}}}$ \hspace{1cm}
    $\vcenter{\hbox{\epsfig{file=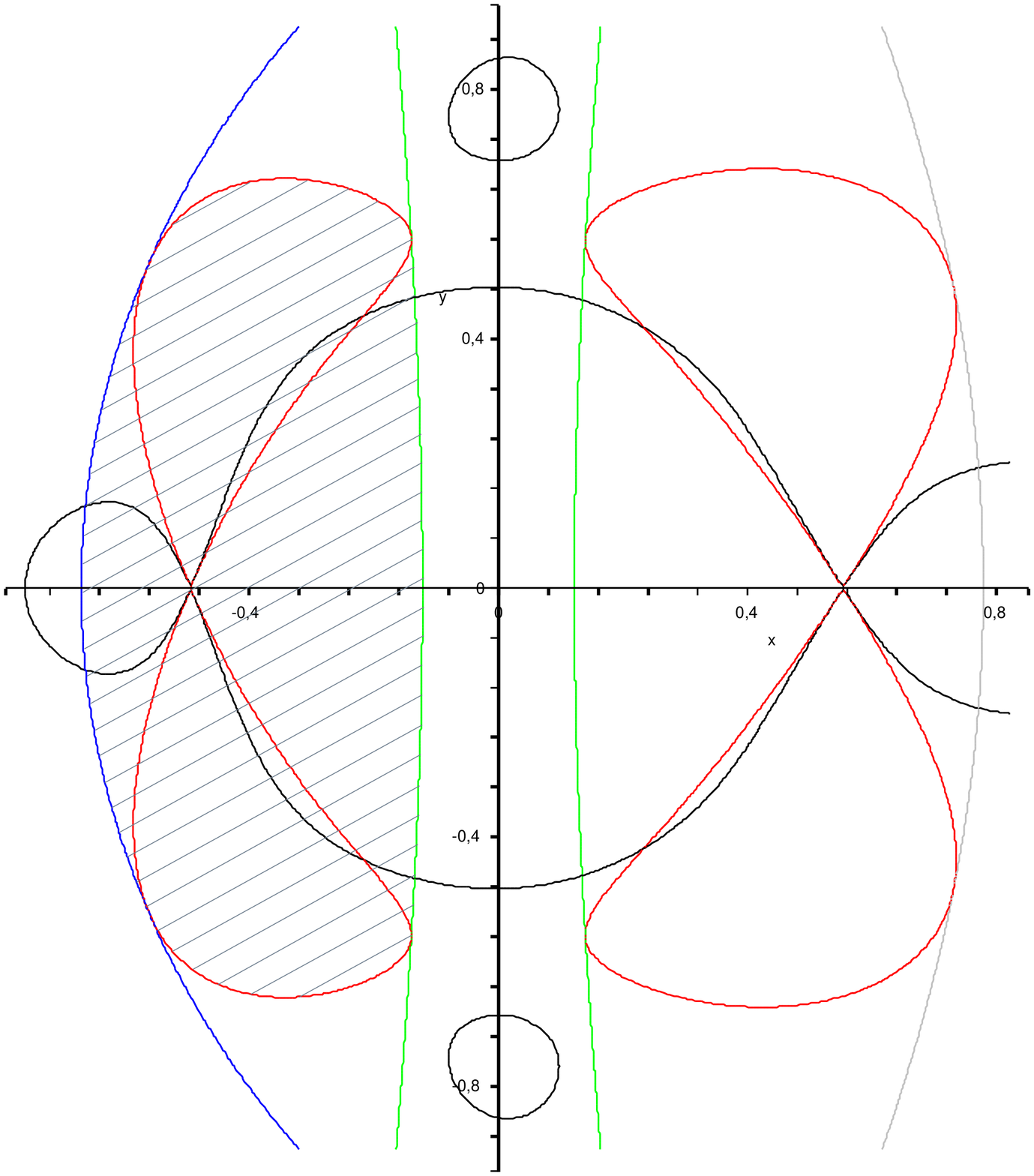,height=6cm,clip=}}}$
   }
  \caption{A configuration of three balls (left) and a planar depiction of a cone of
    directions (right). The direction-sextic is in red, and the Hessian in black.
    The hatched region corresponds to directions of transversals to the
    three balls in the order $B_1 \prec B_0 \prec B_2$, the blue arc corresponding
    to inner special bitangents of balls $B_0$ and $B_2$ and the green arc
    corresponding to inner special bitangents of balls $B_0$ and $B_1$.}
  \label{fig:hatch}
\end{figure}

\paragraph{Cone of directions.}
The \emph{cone of directions} $K(B_0B_1B_2)$ of $B_0,B_1, B_2$ is the set of directions of all oriented line
transversals to these balls which meet them in the stated order: $B_0 \prec B_1 \prec B_2$. The boundary of
$K(B_0B_1B_2)$ consists of \cite[Lemma~9]{CGHP} certain arcs of the direction-sextic $\sigma$, and certain arcs
of directions of \emph{inner special bitangents} i.e. tangents to two of the balls passing through their inner
similitude center \cite{HCV}. Figure~\ref{fig:hatch} offers an illustration of a cone of directions. The plane
of the picture must be conceived as an affine piece $R^2\subset P_2$.

\medskip \noindent
We recall the fact that a common tangent (here called bitangent) for two disjoint spheres (more precisely, the
boundary of two disjoint balls) passes through their inner similitude center if and only if it is contained in a
common tangent plane which has the two spheres on opposite sides. If a transversal for the two balls has the
direction of an inner special bitangent, it must actually be that bitangent. The cone of directions for a pair
of disjoint balls is bounded precisely by their inner special bitangents. In $P_2$ they trace a (circular)
conic.

\medskip \noindent
The points of $\sigma$ that appear on the boundary $\partial K(B_0B_1B_2)$ can be characterized as follows:

\begin{prop}\label{simplex}
The direction of a tritangent $\ell$ meeting the three balls in the prescribed order belongs to $\partial
K(B_0B_1B_2)$ if and only if $\ell$ intersects the triangle of centers $c_0c_1c_2$.
\end{prop}

\medskip \noindent
{\em Proof:}\ The set of directions of common transversals to disjoint balls is a proper subset of $P_2$.

\medskip \noindent
Assume that $\ell$ is neither parallel to the plane of centers, nor contained in it.

\medskip \noindent
If the intercept point of $\ell$ with the plane of centers lies outside the triangle of centers, there exists an
edge which has a center on the other side. This still holds for the projected configuration on $\ell^{\perp}$.
When moving along the perpendicular closer to the projected edge, all distances to projected centers decrease.
This shows that there are lines parallel to $\ell$ ``stabbing'' the open balls, and therefore the direction of
$\ell$ is not on the boundary. On the other hand, when the tritangent $\ell$ intersects the triangle of centers,
and we follow the projection along $\ell$ on $\ell^{\perp}$, there is no distance decreasing motion for all
distances to the (projected) vertices, for this would decrease all areas over edges, while these areas have a
constant sum. Thus the direction of $\ell$ is on $\partial K(B_0B_1B_2)$.

\medskip \noindent
 In other words, the general case follows from the elementary statement
that given a (top dimensional) simplex in a Euclidean space, and a
point, the balls centered at the vertices of the simplex and passing
through the given point intersect only at that point when it belongs
to the simplex, but have an intersection with non-empty interior when
the point is outside the simplex.

\medskip \noindent
If $\ell$ is parallel to the plane of centers (but not contained in it), we may consider any parallel plane
which is closer to $c_0c_1c_2$ than $\ell$ is, and find in this plane transversals to the open balls parallel to
$\ell$. Thus, $\ell$ cannot be on the boundary.

\medskip \noindent
Finally, if $\ell$ is in the plane of centers, and we look at the ``section configuration'' traced in that
plane, {\em either} all three discs are on one side of $\ell$ and then $\ell$ does not cross the triangle of
centers and is not on the boundary, {\em or} $\ell$ has two discs on one side and the third on the other, must
cross the triangle of centers, is actually an inner special bitangent for two pairs of balls (and an outer
special bitangent for the third pair), and is thus necessarily on the boundary.\ \ $\Box$

\begin{prop} \label{closure}
For three disjoint balls, we have:

\medskip
(i)\ \  the cone of directions $K(B_0B_1B_2)$ consists of a single point if and only if there is a tritangent
contained in the plane of centers, and tracing in it a pinned planar configuration i.e. the disc traced by $B_1$
is on the other side of the tritangent than the discs traced by $B_0$ and $B_2$;

\medskip
(ii)\ in all other cases, the cone of directions $K(B_0B_1B_2)$ is the closure of its interior.
\end{prop}

\medskip \noindent
{\em Proof:}\  (i)\ Sufficiency: the plane intersecting the plane of centers along the tritangent, and
perpendicular to it, will have $B_1$ on one side, and $B_0$ and $B_2$ on the other. An oriented transversal
meeting $B_0$ first, then $B_1$, and then $B_2$ must be contained in this separating perpendicular plane, and
thus coincide with the given tritangent. Necessity is covered by our arguments in (ii).

\medskip
(ii)\ Suppose we are not in case (i), and the centers are not aligned. If we have a transversal $\ell$ with
direction belonging to the boundary of $K(B_0B_1B_2)$, we may assume the transversal is not in the plane of
centers, since a non-pinned planar case is clear. But then $\ell$ and its reflection in the plane of centers
define a plane perpendicular to the latter, and all lines between them (passing through their intersection) have
directions belonging to the interior, because all distances from centers decrease.

\medskip \noindent
The case of collinear centers is trivial: there is only one geometric permutation (given by the line of centers)
and the cone of directions is a disc-like region bounded by a conic.\ \ $\Box$

\medskip \noindent
{\bf Corollary of proof:}\ Cone of directions and connected components of transversals for three disjoint balls
in $R^3$  are {\em contractible}.

\medskip \noindent
Indeed, the argument above shows that we may contract first to the segment in $K(B_0B_1B_2)$ consisting of
directions in the plane of centers, and then contract this segment.

\medskip \noindent
Obviously the same holds true at the level of the connected components in the space of transversals.\ \ $\Box$

\paragraph{Hessian and flexes.} The \emph{Hessian} of $\sigma$ is
defined as the determinant of the matrix of second derivatives:

$$    H(\sigma)=H(\sigma)(u)=\det{\left( \frac{\partial^2\sigma}{\partial u_i
      \partial u_j}\right)}$$

\noindent The Hessian curve (or simply Hessian) is the projective curve defined by the zero-set of this
determinant.

\medskip \noindent
The Hessian of a direction-sextic for three balls in $R^3$ is thus an algebraic curve of degree twelve. The
intersection between $\sigma$ and its Hessian $H(\sigma)$ consists of all singular points of $\sigma$ and all
{\em flexes} of $\sigma$ \cite{BK}.

\section{Outline of the proof}\label{sec:outline}

\medskip \noindent
For $d=2$ the convexity theorem is elementary, and for $d\geq 3$ it is easily reduced to the case of three
disjoint balls in $R^3$. The {\em key property} used to settle this case is the following:

\begin{prop}\label{flexes}
For disjoint balls $B_0,B_1,B_2$, any arc of their direction-sextic $\sigma$ which belongs to the boundary
$\partial K(B_0B_1B_2)$ contains no flex or singularity of $\sigma$ between its endpoints.
\end{prop}

\medskip \noindent
The convexity of the cone of directions $K(B_0B_1B_2)$ can then be inferred from the known fact that a simple
$C^1$-loop in $R^2\subset P_2$ with no inflection (in Euclidean terms: with positive curvature on its algebraic
arcs) bounds a convex interior~\cite{Top}.

\medskip \noindent
Thus, what is essential for this approach, is to obtain sufficient control over the flexes of $\sigma$. At first
sight, the fact that the intersection of $\sigma$ and the Hessian $H(\sigma)$ in $P_2(C)$ has (counting
multiplicities) $6\times 12=72$ points leaves little hope for the possibility of ``tracking'' all flexes.
However, there is another way to exploit the Hessian: fix a direction and consider the ball configurations which
have a tritangent with that direction and give the same planar configuration of four points when projecting,
tangent and centers, on some orthogonal plane; express the Hessians of the corresponding direction-sextics and
then ask which may vanish for the given direction.

The important point is that one can anticipate, from the form of the equations, that the computations must
result in polynomials of low degree, which will be subject, in their turn, to geometrical control.

The unfolding of this scenario is presented in the next section and involves a certain amount of explicit
computations. Although no part is too complicated to be done by hand, we have relied on Maple \cite{M} in a
few instances, which are documented in the Appendix.

\section{Details of the proof}\label{sec:computation}

\subsection{Probing for flexes}

Following Proposition~\ref{simplex}, we need only consider directions of tangents to the three balls that cross
the triangle of centers and are not directions of inner special bitangents.  When projecting along such a
tangent on a perpendicular plane, the projected centers form a triangle containing the point image of the
tangent as an interior point. One may start with the latter planar configuration, a triangle and an interior
point, and ask: what ball configurations yield this picture (by projection along a common tangent intersecting
at the interior point)? Since the radii of the balls are given, one has only to ``lift'' the vertices of the
triangle in the normal direction and obtain all the asked for configurations.

\medskip \noindent
We equip $R^3$ with a frame such that the triangle lies in the plane $e_3^{\perp}\subset R^3$ and has its
vertices at $\tilde{c}_0=0, \tilde{c}_1, \tilde{c}_2$, with the understanding that there is a point inside, with
squared distances $s_i$ to these vertices. Then, we use three real parameters, $x_0, x_1$ and $x_2$, to describe
the possible positions of the three centers:

$$   c_0=\tilde{c}_0+x_0 e_3,\ \ c_1=\tilde{c}_1+x_1e_3, \ \ c_2=\tilde{c}_2+x_2e_3 $$

\medskip \noindent
We use Proposition~\ref{sextic} to express the corresponding direction-sextic $\sigma$ and its Hessian
$H(\sigma)$ as functions of $x=(x_0,x_1,x_2)\in R^3$ depending on
$\tilde{c}_0,\tilde{c}_1,\tilde{c}_2,s_0,s_1,s_2$. Proposition~\ref{flexes} is now equivalent to proving that

$$ H(\sigma)(0,0,1)\neq 0 $$

\noindent holds for all initial data (triangle and interior point) and all $(x_0,x_1,x_2)$ corresponding to
disjoint balls.

\subsection{A quadric and a quartic}

\medskip \noindent
We have reduced the probe for flexes to the study of a polynomial function of $x$ (and parameters) which can be
explicitly computed. For the Maple procedure we used, see Appendix.

\medskip \noindent
The parameters involved are the following:

$$   \tilde{c}_0=(0,0,0),\ \tilde{c}_1=(a,0,0),\ \tilde{c}_2=(b,c,0) $$

\noindent the triangle of centers $(\tilde{c}_0,\tilde{c}_1,\tilde{c}_2)$ having interior point:

$$   p = \frac{\sum p_i\tilde{c}_i}{\sum p_i} = \frac{p_1\tilde{c}_1+p_2\tilde{c}_2}{\sum p_i},\ \ p_0, p_1, p_2
> 0 $$

\medskip \noindent
Let $v_k=p-\tilde{c}_k$. Then \   $s_k = r_k^2 = <v_k,v_k>$.

\medskip
\noindent The computation gives the result:

$$   H(\sigma)(0,0,1) = \frac{2^{12}5^2a^6c^6}{(\sum p_i)^5} [H_2(x)+H_4(x)] $$

\noindent where $H_2$ and $H_4$ have degree respectively 2 and 4 in $x=(x_0,x_1,x_2)$:

$$   H_2=H_2(x)=-a^2c^2\ \prod p_k \sum p_ip_j(x_i-x_j)^2 $$

$$   H_4=H_4(x)= \sum p_k^3s_k(x_i-x_k)^2(x_j-x_k)^2 $$

\noindent with cyclic notation for $\{i,j,k\}=\{0,1,2\}$. Thus, away from $(0,0,0)$, $H_2$ is negative and $H_4$
is positive. The aim is now to show that the assumption of disjoint balls is enough to ensure the positivity of
$H_2+H_4$.

\subsection{Hyperboloid and octant}

\medskip \noindent
We can further transform these expressions by retaining as parameters the (positive numbers) $p_i$ and
$q_j=p_jr_j$, and renaming the squares $z_k=(x_i-x_j)^2$. This gives:

$$   H_2=H_2(z)=-a^2c^2\prod p_k \sum p_ip_jz_k $$

$$   H_4=H_4(z)= \sum p_kq_k^2z_iz_j $$

\noindent {\bf From now on, assume that} $\sum p_i = 1$. We have to replace $\Delta=a^2c^2$, which is four times
the squared area of the triangle $\tilde{c}_0,\tilde{c}_1,\tilde{c}_2$, by its expression in terms of $p_i$ and
$q_j$.

\begin{lemma}
  We have:

$$  \Delta=a^2c^2=\frac{Q}{4\prod p_k^2}, \ \  \mbox{with} \ Q = \sum (2q_i^2q_j^2-q_k^4) $$
\end{lemma}

\medskip \noindent {\em Proof:}\  This is an elementary computation, which may be conducted as follows. By the
  definition of $v_i$, we have

  $$   \sum p_iv_i=0 $$

\noindent
  From $<\sum p_iv_i,v_j>=0$, we obtain a linear system for $<v_i,v_j>,\ i\neq j$:

$$    p_i<v_i,v_k>+p_j<v_j,v_k> = -p_k <v_k,v_k>=-p_ks_k $$

\noindent
  with solutions:

  $$   <v_i,v_j>=\frac{p_k^2s_k-p_i^2s_i-p_j^2s_j}{2p_ip_j}=\frac{q_k^2-q_i^2-q_j^2}{2p_ip_j} $$

  \medskip \noindent
  Four times the squared area of a triangle $p,\tilde{c}_i,\tilde{c}_j$ is a Gram determinant:

  $$  \left| \begin{array}{cc}
       <v_i,v_i> & <v_i,v_j>\\
       <v_i,v_j> & <v_j,v_j>
     \end{array} \right| = s_is_j-<v_i,v_j>^2=\frac{Q}{4p_i^2p_j^2} $$

  \noindent
  where $Q=\sum (2q_i^2q_j^2-q_k^4)$. Hence the area of the triangle $\tilde{c}_0,
  \tilde{c}_1,\tilde{c}_2$ is:

  $$   \frac{1}{4}Q^{1/2}\sum \frac{1}{p_ip_j}=\frac{Q^{1/2}}{4\prod p_k} $$

  \noindent
  resulting in:

  $$   \Delta=a^2c^2=\frac{Q}{4\prod p_k^2} $$ $\Box$

\medskip \noindent
Several new substitutions will be in order for the study of $H_2+H_4$. Since a positive factor won't affect sign
considerations, we'll use the {\bf symbol} $*H$ for any positive multiple of $H_2+H_4$. We have found above:

$$   *H=*H(z)=-\frac{1}{4} Q \sum \frac{z_k}{p_k} + \sum p_kq_k^2z_iz_j $$

\noindent with the shorthand $Q=\sum (2q_i^2q_j^2-q_k^4)$. We put $p_ip_jz_k=q_k^2w_k$, and obtain (up to a
positive factor):

$$   *H=*H(w)=-\frac{1}{4}Q \sum q_k^2w_k + \prod q_k^2 \sum w_iw_j $$

\medskip \noindent
With one more positive rescaling, and $a_k=\frac{Q}{4q_i^2q_j^2}$, we have:

$$   *H=*H(w)=\sum w_iw_j-\sum a_kw_k $$

\medskip \noindent
We can turn now to the conditions expressing the fact that the spheres with centers $c_i=\tilde{c}_i+x_ie_3$ and
radii $r_i$ are disjoint. They are:

$$   z_k=(x_i-x_j)^2 > (r_i+r_j)^2 - \delta_{ij}=(r_i+r_j)^2-<v_i-v_j,v_i-v_j> $$

\noindent that is:

$$   z_k > \frac{q_k^2-(q_i-q_j)^2}{p_ip_j} $$

\medskip \noindent
In $w$-coordinates, the {\bf ``disjointness conditions''} become

$$   w_k > 1-\left(\frac{q_i-q_j}{q_k}\right)^2 $$

\medskip \noindent
Note that from $\sum p_iv_i=0$ it follows that $q_k=\|p_iv_i\| > 0$ are the edges of a triangle, and therefore
the latter expressions are positive by the triangle inequality.

\medskip \noindent
The purpose now is to study the position of the octant defined by the ``disjointness conditions'' relative to
the affine quadric in $R^3$ defined by $*H(w)=0$. We use first a translation by $\beta$, in order to absorb the
linear part in $*H$:

$$   *H=*H(w)=\sum (w_i-\beta_i)(w_j-\beta_j)-\sum \beta_i\beta_j $$

\noindent requesting:

$$   \beta_i+\beta_j= a_k,\ \mbox{that is}\ \ \beta_k=\frac{1}{2}(a_i+a_j-a_k) $$

\medskip \noindent
This makes

$$   \sum \beta_i\beta_j=\frac{1}{4}\sum (a_k+a_i-a_j)(a_k-a_i+a_j)=
     \frac{1}{4}\sum (2a_ia_j-a_k^2) $$

\noindent and results in

$$   \sum \beta_i\beta_j=\frac{1}{4} \left(\frac{Q}{4\prod q_k^2}\right)^2 \sum
     (2q_i^2q_j^2-q_k^4)= \frac{Q^3}{4^3\prod q_k^4} > 0 $$

\medskip \noindent
Thus, with translated coordinates $t_k=w_k-\beta_k$ we have a \textbf{hyperboloid with two sheets}:

$$   *H=*H(t)=\sum t_it_j-\frac{Q^3}{4^3\prod q_k^4}=0 $$

\noindent which lies on the positive side of its asymptotic cone $\sum t_it_j=0$.

\begin{lemma}
  $\sum t_it_j=0$ is a circular cone with axis $t_0=t_1=t_2$. The two components of its smooth points
  circumscribe the positive and negative open octants, which are both contained in the
  positive part $\sum t_it_j > 0$.
\end{lemma}

\medskip \noindent
The open octant defined by our ``disjointness conditions'' $w_k >
1-(\frac{q_i-q_j}{q_k})^2$ is a translate of the open positive octant, and its
position relative to the hyperboloid $*H(w)=0$ is determined by the position of
its vertex $V$. Continuing to refer here to $w$-coordinates, we have:

\begin{lemma}
  The point $V=(1-(\frac{q_i-q_j}{q_k})^2)_{0\leq k\leq 2}$ is on the ``positive
  side'' of the hyperboloid $*H(w)=0$ and on the ``positive side'' of the plane $\sum
  t_k=\sum (w_k-\beta_k)=0$, that is:

  $$   *H(V) > 0 \ \ \ \mbox{and} \ \ \ \sum \left(1-\left(\frac{q_i-q_j}{q_k}\right)^2\right) >
       \frac{Q}{8\prod q_k^2}\sum q_k^2 $$
\end{lemma}

\medskip \noindent
{\em Proof:} \ A Maple assisted computation (see Appendix) shows that
  $*H(V)$ factors as

  $$   *H(V)=\frac{3\prod (q_i+q_j-q_k)^2}{4\prod q_k^2} $$

  \noindent
  from which the first inequality follows.

\medskip \noindent
  The second inequality, which determines
  on which of the two components of the ``positive side'' of the hyperboloid $V$
  lies, is satisfied for $q_0=q_1=q_2$, and by continuity, must be satisfied for
  any other triangle edges, since vertex $V$ cannot ``jump'' from one component to
  the other. \ \ $\Box$

\medskip \noindent
It is now clear, geometrically, that the octant where the ``disjointness conditions'' are satisfied and the
hyperboloid indicating a flex or a singularity for the corresponding configuration {\em have no point in
common}. This completes the proof of Proposition~\ref{flexes}.

\subsection{Convexity of the cone of directions}

\medskip \noindent
We consider now three {\bf disjoint} closed balls $B_0,B_1,B_2$ described by parameters: centers $c_0,c_1,c_2$
and radii $r_0,r_1,r_2$. We shall prove first the convexity of any cone of directions in the {\em generic} case
i.e. when the centers are in the complement of a proper algebraic subset. Then, we'll show that the generic case
implies the general case.

\begin{lemma}
The direction cone $K(B_0B_1B_2)$ of a generic triple of disjoint balls in $R^3$ is strictly convex.
\end{lemma}

\medskip \noindent
{\em Proof:}\ Genericity allows us to assume that the direction-sextic $\sigma$ is non-singular at all its
contacts with any of the three conics determined by inner special tangents. Then, these contacts are tangency
points, and if we start at some point of, say $\partial K(B_0B_1B_2)$, and follow the boundary curve, we obtain,
by Proposition~\ref{flexes}, a differentiable simple loop of class $C^1$, which is, locally, always on the same
side of its tangent. For any affine plane $R^2\subset P_2$ covering the loop, and any Euclidean metric in it,
this means positive curvature on all its algebraic arcs, and this implies \cite{Top} the fact that our simple
loop bounds a compact convex set. In fact {\em strictly convex}, because of non-vanishing curvature. By
Proposition~\ref{closure} and its Corollary, this strictly convex set is $K(B_0B_1B_2)$. \ \ $\Box$

\medskip \noindent
The passage from the generic case to the general case is based on:

\begin{lemma}
Let ${\cal B}=(B_0,B_1,B_2)$ be a configuration of three disjoint closed balls, and suppose $K(B_0B_1B_2)$ has
non-empty interior. If ${\cal B}$ is the limit of a sequence of configurations ${\cal B}^{(\nu)}$ with a convex
corresponding cone of directions, then $K(B_0B_1B_2)$ is convex as well.
\end{lemma}

\medskip \noindent
{\em Proof:}\ By Proposition~\ref{closure}, it is enough to prove that, for any two points in the interior, the
(geodesic) segment joining them is contained in $K(B_0B_1B_2)$.

\medskip \noindent
Take two interior points. By assumption, for sufficiently large $\nu$, the segment joining them is contained in
all corresponding cones for ${\cal B}^{(\nu)}$. Consider one point of the segment, and project the sphere
configuration along the direction defined by the point, on a perpendicular plane. We have to prove that the
disks representing the projected balls have at least one point in common.

Suppose they don't. Then so would discs with the same centers and radii increased by a small $\epsilon > 0$. But
then we can find, for sufficiently large $\nu$, configurations ${\cal B}^{(\nu)}$ with centers projecting less
than $\epsilon/2$ away from those of ${\cal B}$, and corresponding radii with less than $\epsilon/2$
augmentation. Then the point of the segment cannot be in the respective cones of directions: a
contradiction. \ \ \ $\Box$

\medskip \noindent
The convexity result generalizes to arbitrary $n$ and $d$  as follows:

\medskip \noindent
{\bf Proof of Theorem~\ref{thm:convexity}}\  Recall that, for any collection of balls in $R^3$, a direction will
be realized by some transversal if and only if the orthogonal projection of the balls on a perpendicular plane
has non-empty intersection. By Helly's Theorem in the plane, the direction cone for a sequence of $n \geq 3$
balls is the intersection of the direction cones of all its triples. Thus, the direction cone of $n$ ordered
$3$-dimensional disjoint balls is strictly convex for any $n$.

\medskip \noindent
Given a sequence ${\cal S}$ of $n$ disjoint balls in $R^d$, let $K$ be its direction cone for a prescribed order
of intersection. Let $u$ and $v$ be two directions in $K$, $\ell_u$ and $\ell_v$ be two corresponding line
transversals and let $E$ denote the $3$-dimensional affine space these two lines span (or a $3$-space containing
their planar span, should the lines be coplanar).

\medskip \noindent
$E \cap {\cal S}$ is a collection of $3$-dimensional disjoint balls whose corresponding direction cone is convex
on $S^2$. Thus, for any direction on the small arc of great circle joining $u$ and $v$ there exists an
order-respecting transversal to ${\cal S}$, because it already exists in $E$. It follows that $K$ is convex, and
again, from the three dimensional case, strictly convex. \ \ $\Box$

\begin{figure}[!]
  \vspace{-1cm}
  \centerline{a. \epsfig{file=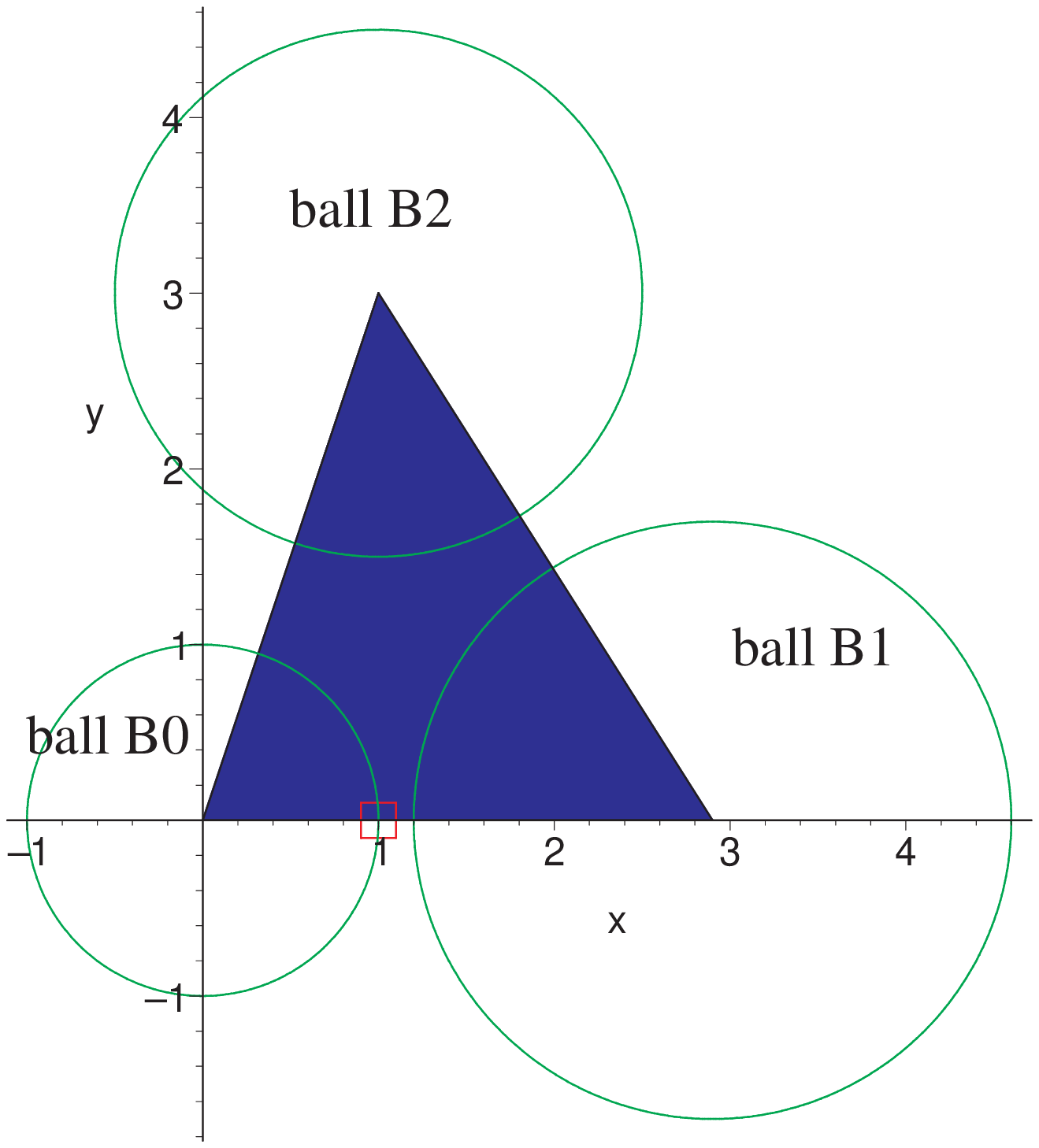,height=6.6cm}}

  \vspace{0.8cm}

  \centerline{$\vcenter{\hbox{b.\epsfig{file=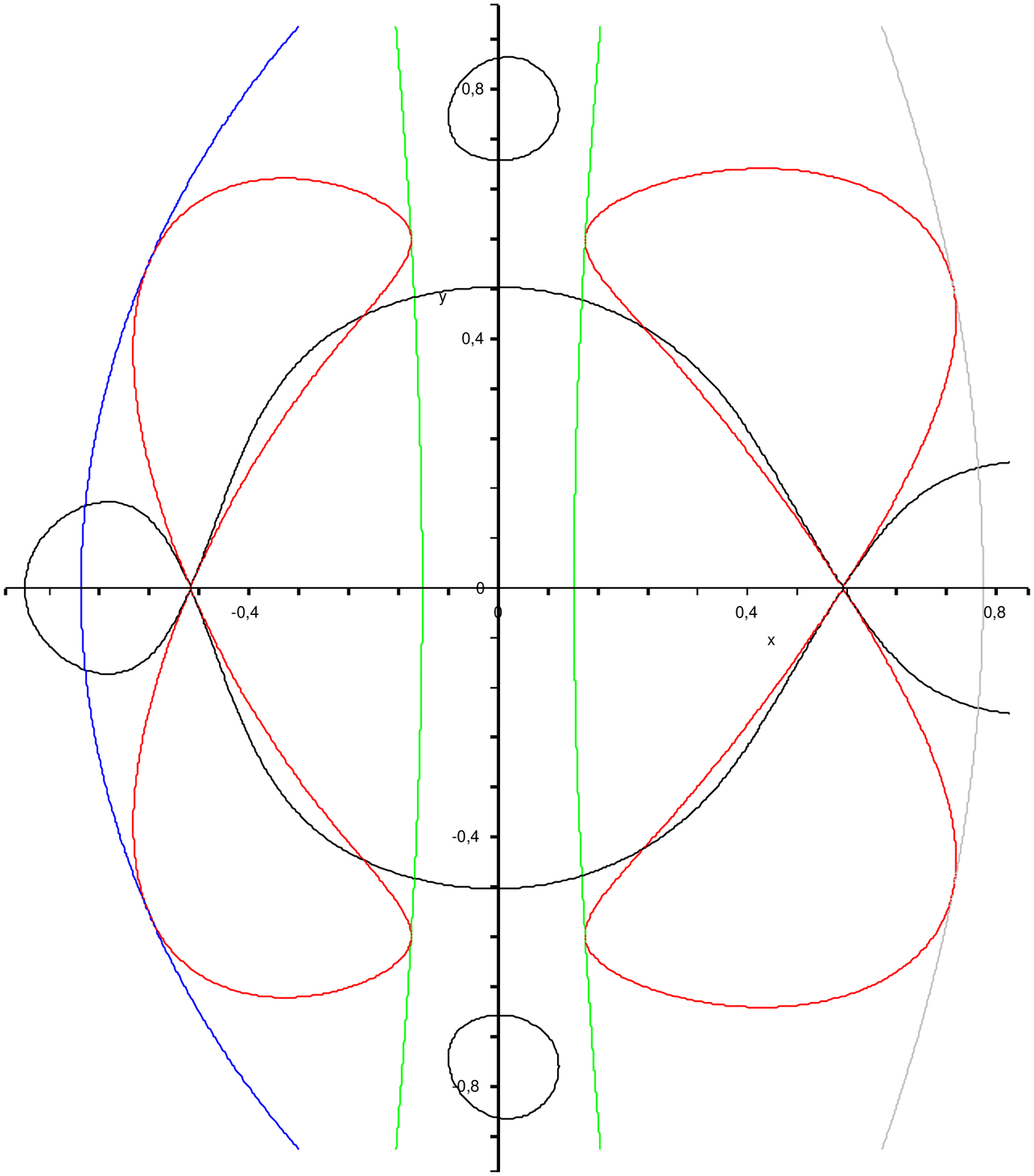,width=5cm,clip=}}}$ \hspace{2mm}
    $\vcenter{\hbox{\epsfig{file=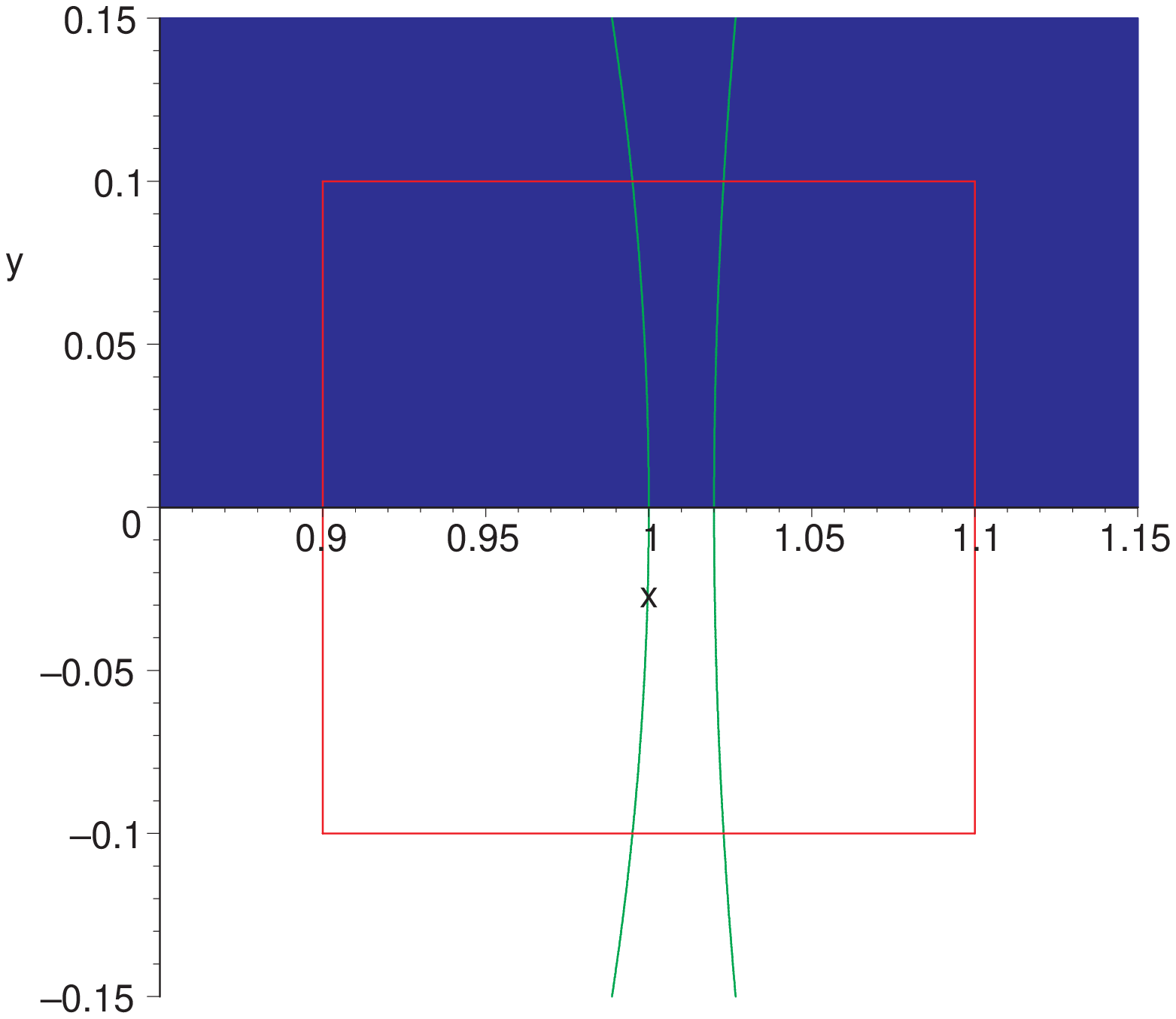,height=2.6cm}}}$ \hspace{0.3cm}
    $\vcenter{\hbox{c.\epsfig{file=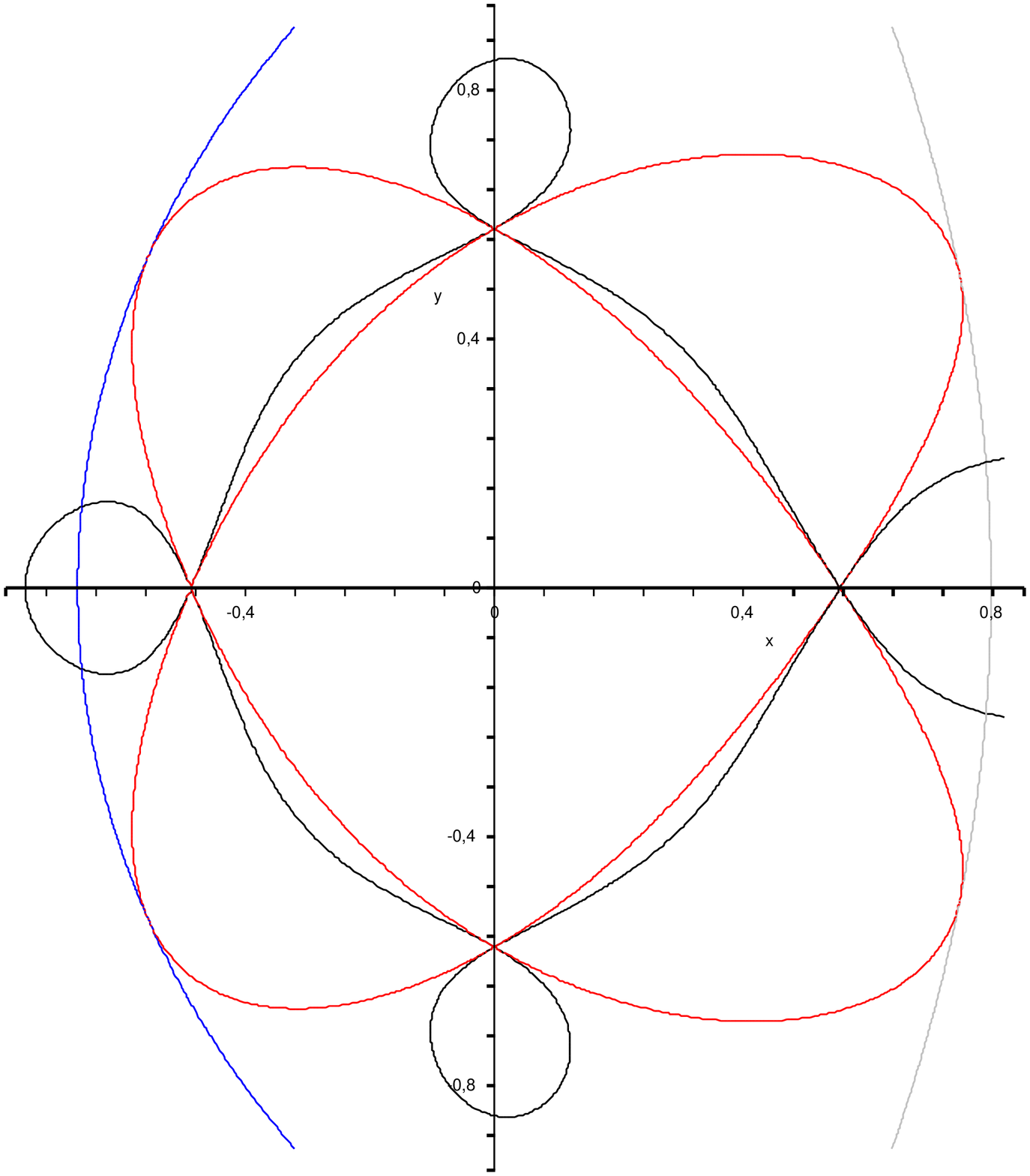,width=5cm,clip=}}}$ \hspace{2mm}
    $\vcenter{\hbox{\epsfig{file=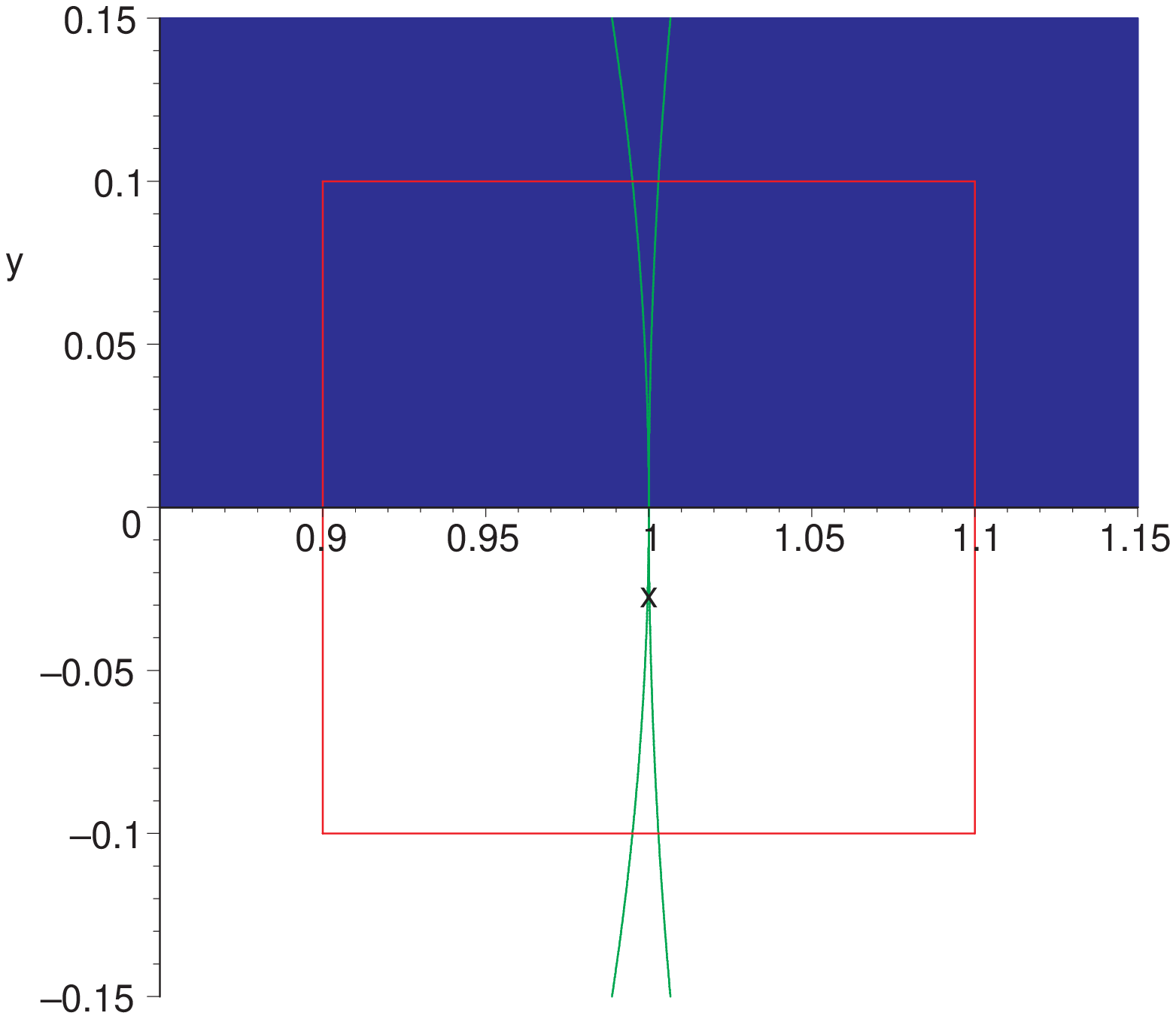,height=2.6cm}}}$}

  \vspace{0.8cm}

  \centerline{$\vcenter{\hbox{d.\epsfig{file=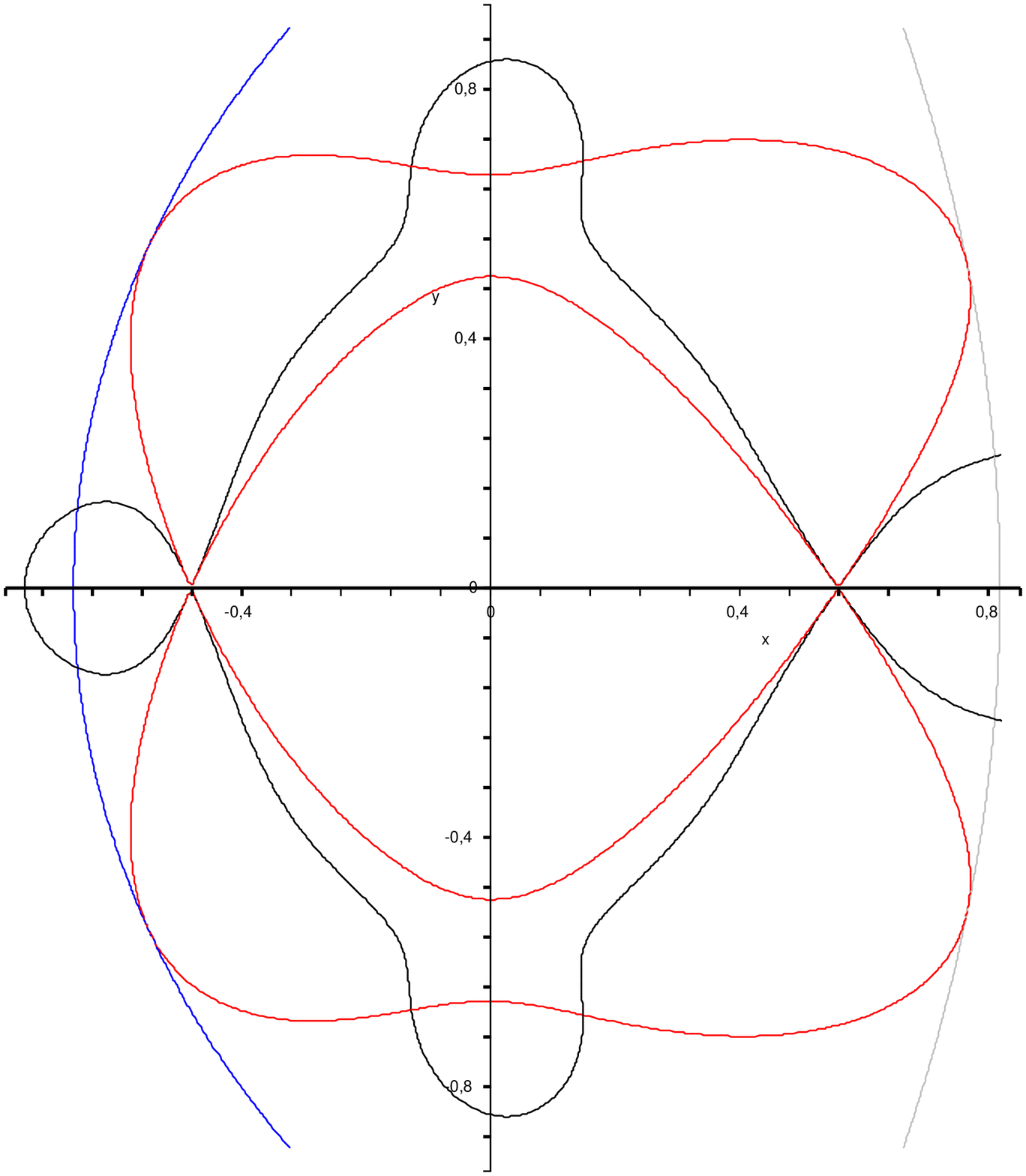,width=5cm,clip=}}}$ \hspace{2mm}
    $\vcenter{\hbox{\epsfig{file=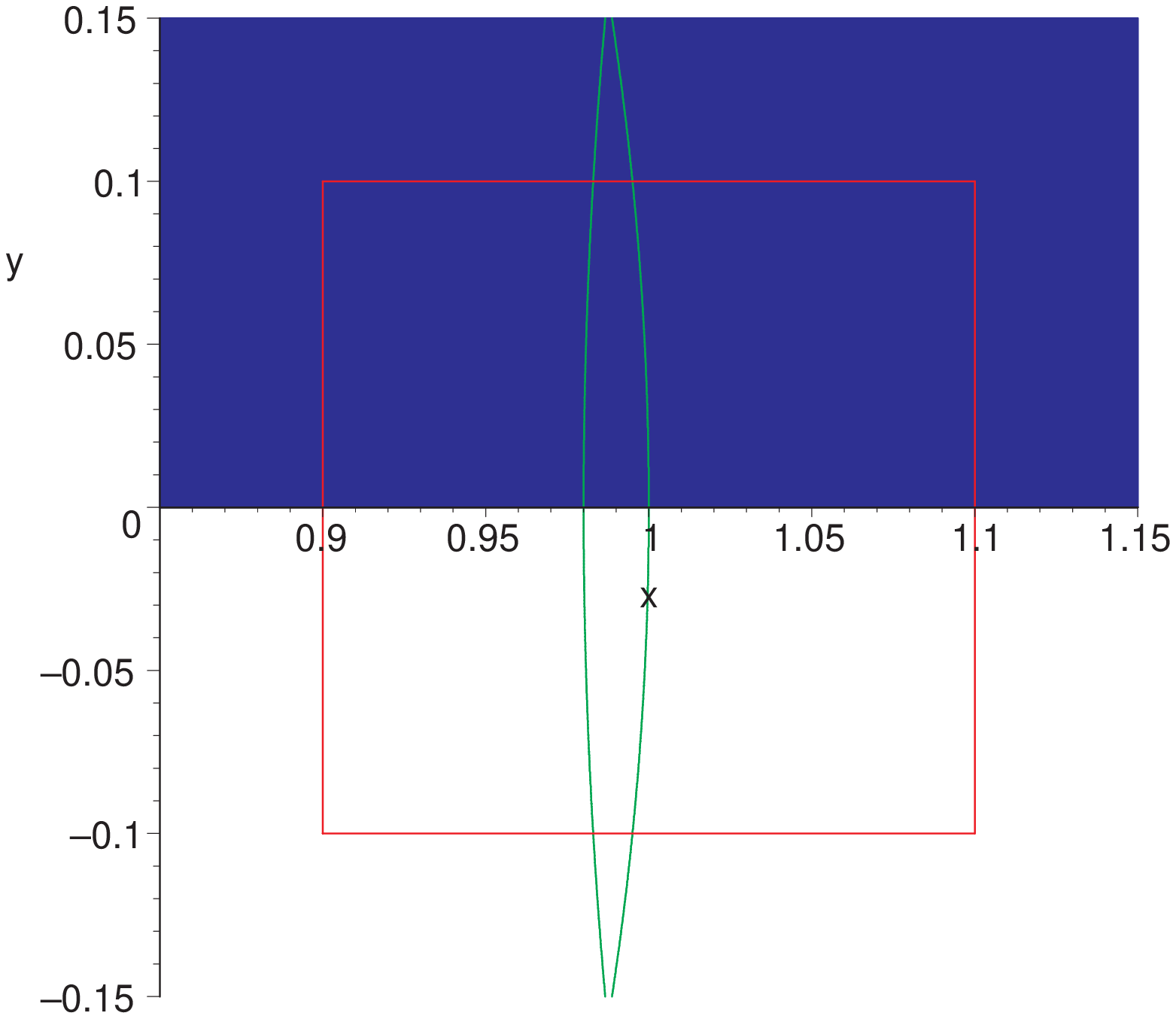,height=2.6cm}}}$}
  \vspace{0.8cm}
  \caption{a. The trace of three disjoint balls on the plane of centers, with ball
    $B_1$ moving on the horizontal axis towards ball $B_0$. The red square is used
    for close-ups below. b. c. d. The direction-sextic (in red), its Hessian (in
    black) and the inner special bitangent conics (in blue, green and gray), when
    balls $B_0$ and $B_1$ are disjoint (b), tangent (c) and intersecting (d).}
  \label{fig:transition}
\end{figure}

\medskip \noindent
Let us emphasize the importance of the assumption that the balls are disjoint. Figure~\ref{fig:transition}
illustrates a transition from convex to non-convex direction cones as three disjoint balls move and allow an
overlap.

\section{Implications}\label{sec:implications}

\medskip \noindent
This section explores some consequences of Theorem~\ref{thm:convexity}. Similar results were proven for the case
of unit balls in \cite{CGHP} and, with Theorem~\ref{thm:convexity}, the proofs carry through. We thus omit all
arguments here and point to the relevant lemmata in \cite{CGHP}.

\subsection{Isotopy and geometric permutations}

\medskip \noindent
An immediate corollary of Theorem~\ref{thm:convexity} is the correspondence of isotopy and geometric
permutations for line transversals to disjoint balls:

\begin{cor}\label{cor:topology}
The set of line transversals to $n$ disjoint balls in $R^d$ realizing the same geometric permutation is
contractible.
\end{cor}

\medskip
\noindent The proof given by Cheong et al. \cite[Lemma~14]{CGHP} for disjoint unit balls immediately extends,
with Theorem~\ref{thm:convexity}, to the case of disjoint balls. Smorodinsky et al. \cite{SMS} showed that in
the worst case $n$ disjoint balls in $R^d$ admit $\Theta(n^{d-1})$ geometric permutations. The same bound thus
applies for the number of connected components of line transversals, improving on the previous bounds of
$O(n^{3+\epsilon})$ for $d=3$ and of $O(n^{2d-2})$ for $d \geq 4$ due to Koltun and Sharir \cite{KoS}. If the
radii of the balls are in some interval $[1, \gamma]$ where $\gamma$ is independent of $n$ and $d$, then the
number of components of transversals is $O(\gamma^{\log \gamma})$, following the bound on the number of
geometric permutations obtained by Zhou and Suri~\cite{ZS}. These results are summarized as follows:

\begin{cor}
In the worst case, $n$ disjoint balls in $R^d$ have $\Theta(n^{d-1})$ connected components of line transversals.
If the radii of the balls are in the interval $[1, \gamma]$, where $\gamma$ is independent of $n$ and $d$, this
number becomes $O(\gamma^{\log \gamma})$.
\end{cor}

\subsection{Minimal pinning configurations}

\medskip \noindent
A \emph{minimal pinning configuration} is a collection of objects having an isolated line transversal that
ceases to be isolated if any of the objects is discarded. An important step in the proof of Hadwiger's
transversal theorem \cite{Had1} is the observation that, in the plane, any minimal pinning configuration
consisting of disjoint convex objects has cardinality $3$. Cheong et al. \cite[Proposition~13]{CGHP} proved that
any minimal pinning configuration consisting of disjoint unit balls in $R^d$ has cardinality at most $2d-1$.
With Theorem~\ref{thm:convexity}, the same holds for disjoint balls of arbitrary radii:

\begin{cor}\label{cor:pinning}
Any minimal pinning configuration consisting of disjoint balls in $R^d$ has cardinality at most $2d-1$.
\end{cor}

\subsection{A Hadwiger-type result}

\medskip \noindent
The ``pure'' generalizations \cite{CGHP,HKL} of Helly's theorem, i.e. without additional constraints on
the ordering \emph{\`a la} Hadwiger, use two ingredients: the convexity of the cone of directions and the fact
that $n \geq 9$ disjoint unit balls have at most $2$ geometric permutations \cite{CGN}. Since the latter is not
true for balls of arbitrary radii \cite{SMS}, such theorems do not generalize immediately to non-unit balls.
Yet, an intermediate result of independent interest in the flavor of Hadwiger's transversal theorem does
generalize:

\begin{cor}\label{cor:Hadwiger}
A sequence of $n$ disjoint balls in $R^d$ has a line transversal if any subsequence of size at most $2d$ has an
order-respecting line transversal.
\end{cor}

\newpage

\section*{Appendix:\ Maple code\label{app:maple}}

\noindent Maple 10 code \cite{M} for computations used in Section 4.

\medskip

\begin{verbatim}

with(LinearAlgebra):
with(VectorCalculus):

u := Vector([u1,u2,u3]):
q := DotProduct(u,u):

###### The vertices of the triangle
c0t := Vector([0,0,0]): c1t := Vector([a,0,0]): c2t := Vector([b,c,0]):

###### The centers of the spheres
c0 := c0t+x0*Vector([0,0,1]): c1 := c1t+x1*Vector([0,0,1]):
c2 := c2t+x2*Vector([0,0,1]):

###### Additional variables
e01 := c1-c0: e02 := c2-c0: e12 := c2-c1:
cp01 := CrossProduct(e01,u): cp02 := CrossProduct(e02,u): cp12 := CrossProduct(e12,u):

###### Coefficients of the matrix defining the direction sextic
t01 := DotProduct(cp01,cp01): t02 := DotProduct(cp02,cp02):
t12 := DotProduct(cp12,cp12):

###### The point in the triangle
po := (p0*c0t+p1*c1t+p2*c2t)/(p0+p1+p2):

###### Squared distances between vertices of the triangle and point inside
s0 := DotProduct(po-c0t,po-c0t): s1 := DotProduct(po-c1t,po-c1t):
s2 := DotProduct(po-c2t,po-c2t):

###### Matrix defining the direction sextic
sigm := Matrix(5,5,[[0,1,1,1,1],[1,0,q*s0,q*s1,q*s2],[1,q*s0,0,t01,t02],
                    [1,q*s1,t01,0,t12],[1,q*s2,t02,t12,0]]):

###### Equation of the sextic
sig := Determinant(sigm):

###### Hessian matrix
Hm := Hessian(sig,[u1,u2,u3]):
Hb := subs(u1=0,u2=0,u3=1,Hm):

###### Hessian curve
H := Determinant(Hb):

###### Divide by positive constant
H := numer(factor(H/(2^(12)*5^2*a^6*c^6))):

###### Decomposition H = H2+H4
H2 := factor(coeff(coeff(coeff(H,x0,2),x1,0),x2,0)*x0^2+
             coeff(coeff(coeff(H,x1,2),x2,0),x0,0)*x1^2+
             coeff(coeff(coeff(H,x2,2),x0,0),x1,0)*x2^2+
             coeff(coeff(coeff(H,x0,1),x1,1),x2,0)*x0*x1+
             coeff(coeff(coeff(H,x0,1),x2,1),x1,0)*x0*x2+
             coeff(coeff(coeff(H,x1,1),x2,1),x0,0)*x1*x2):

H4 := factor(expand(H-H2)):

###### Substitute for yi, with yi = xj-xk
H2y := factor(p0*p1*p2*(algsubs(x1-x2=y0,factor(coeff(H2/(p0*p1*p2),p0,0)))+
                        algsubs(x2-x0=y1,factor(coeff(H2/(p0*p1*p2),p1,0)))+
                        algsubs(x0-x1=y2,factor(coeff(H2/(p0*p1*p2),p2,0))))):

H4y := algsubs(numer(s0t)=r0^2,algsubs(x2-x0=y1,algsubs(x0-x1=y2,
                             factor(coeff(H4,p0,3)),exact),exact),exact)*p0^3+
       algsubs(numer(s1t)=r1^2,algsubs(x1-x2=y0,algsubs(x0-x1=y2,
                             factor(coeff(H4,p1,3)),exact),exact),exact)*p1^3+
       algsubs(numer(s2t)=r2^2,algsubs(x1-x2=y0,algsubs(x2-x0=y1,
                             factor(coeff(H4,p2,3)),exact),exact),exact)*p2^3:

###### Substitute qk = pk*rk and zk = yk^2
H2z := algsubs(y0^2=z0,algsubs(y1^2=z1,algsubs(y2^2=z2,H2y))):
H2z := algsubs(p0*r0=q0,algsubs(p1*r1=q1,algsubs(p2*r2=q2,H2z))):

H4z := algsubs(y0^2=z0,algsubs(y1^2=z1,algsubs(y2^2=z2,H4y))):
H4z := algsubs(p0*r0=q0,algsubs(p1*r1=q1,algsubs(p2*r2=q2,H4z))):

###### Substitute a^2 c^2 = Q/(4*p0^2*p1^2*p2^2)
H2z := algsubs(a^2*c^2=Q/(4*p0^2*p1^2*p2^2),H2z):

###### Substitute pi*pj*zk = qk^2*wk and multiply by a positive scalar
H2w := subs(z2=q2^2*w2/p0/p1,subs(z1=q1^2*w1/p0/p2,subs(z0=q0^2*w0/p1/p2,H2z))):

H4w := subs(z2=q2^2*w2/p0/p1,subs(z1=q1^2*w1/p0/p2,subs(z0=q0^2*w0/p1/p2,H4z))):

H2w := factor(4*p0*p1*p2*H2w):

H4w := factor(4*p0*p1*p2*H4w):

###### The vertex V of the disjointness conditions
v0 := 1-(q1-q2)^2/q0^2: v1 := 1-(q2-q0)^2/q1^2: v2 := 1-(q0-q1)^2/q2^2:

###### Evaluation of H at this vertex
print(factor(subs(w0=v0,w1=v1,w2=v2,
            Q=2*q0^2*q1^2+2*q0^2*q2^2+2*q1^2*q2^2-q0^4-q1^4-q2^4,H2w+H4w))):

###### Evaluation of the plane t0+t1+t2 at the vertex V for q0=q1=q2
plane := w0+w1+w2-Q*(q0^2+q1^2+q2^2)/8/(q0^2*q1^2*q2^2):

print(subs(q1=q0,q2=q0,subs(w0=v0,w1=v1,w2=v2,
            Q=2*q0^2*q1^2+2*q0^2*q2^2+2*q1^2*q2^2-q0^4-q1^4-q2^4,plane))):
\end{verbatim}

\end{document}